\documentclass[12pt,a4paper]{article}
\pdfoutput=1
\usepackage[english]{babel}
\usepackage{verbatim}
\usepackage{bm}

\usepackage{amsmath}
\usepackage{amssymb}
\usepackage{amsfonts}
\usepackage{amsbsy}

\usepackage{graphicx}
\usepackage{graphics}
\usepackage{subfigure}
\usepackage{epsfig}

\usepackage{array}
\usepackage{booktabs}
\usepackage{caption}
\usepackage{tabularx}

\newcommand{\R}{{\mathbb R}}

\newcommand{\Z}{\mathbb Z}

\begin{document}

\title{Semi-Lagrangian methods \\
for parabolic problems in divergence form}

\author{Luca Bonaventura$^{(1)}$,\ \  Roberto Ferretti$^{(2)}$}

\maketitle

\begin{center}
{\small
$^{(1)}$ MOX -- Modelling and Scientific Computing, \\
Dipartimento di Matematica ``F. Brioschi'', Politecnico di Milano \\
Via Bonardi 9, 20133 Milano, Italy\\
{\tt luca.bonaventura@polimi.it}
}
\end{center}

\begin{center}
{\small
$^{(2)}$
Dipartimento di Matematica e Fisica\\
Universit\`a degli Studi Roma Tre\\
 L.go S. Leonardo Murialdo 1,
00146, Roma, Italy\\
{\tt ferretti@mat.uniroma3.it}
}
\end{center}

\date{}


\vspace*{0.5cm}


\vspace*{0.5cm}

\pagebreak

\abstract{Semi-Lagrangian methods have traditionally been developed in the framework of hyperbolic equations, 
but several  extensions of the Semi-Lagrangian approach to diffusion and advection--diffusion problems have been proposed recently. These extensions are mostly based on probabilistic arguments and share the common feature of treating second-order operators in trace form, which makes them unsuitable for  mass conservative models like the classical formulations of turbulent diffusion employed in computational fluid dynamics.
We propose here some basic ideas for treating second-order operators in divergence form. A general framework for constructing consistent schemes in one space dimension is presented, and a specific case of nonconservative discretization is discussed in detail and analysed.  
Finally, an extension to (possibly nonlinear) problems in an arbitrary number of dimensions is proposed. 
Although the resulting discretization approach is only of first order in time, numerical results in a number of test cases 
highlight the advantages of these methods for applications to computational fluid dynamics and their superiority over to more standard
low order time discretization approaches.}

\pagebreak

 \section{Introduction}
 \label{intro} \indent
 
 Semi-Lagrangian (SL) methods are a well established approach in the numerical approximation of  advection dominated problems, that enables to achieve high order accuracy and unconditional stability. They have been
 widely applied in a number of different areas, in particular environmental fluid dynamics and meteorology. A recent detailed presentation of these methods is given in \cite{falcone:2013}, along with a comprehensive literature review. 
 
 Similar strategies, derived from the probabilistic Feynman--Kac representation formula, have been proposed for their extension to parabolic problems, see e.g. \cite{camilli:1995}, \cite{teixeira:1999}, \cite{milstein:2000}, \cite{milstein:2001}, \cite{milstein:2002}, \cite{ferretti:2010}.  Due to their probabilistic interpretation, all these SL methods are derived for parabolic problems written in trace form as
 \begin{equation}
 u_t = \sum_{i,j}\mu_{i,j} u_{x_ix_j}.
 \label{par_noncons}
 \end{equation}
However, in most applications where a variable diffusivity is employed, the corresponding parabolic problem is  formulated in divergence form
\begin{equation}
 u_t = \sum_{i}\Big(\sum_j\nu_{i,j} u_{x_j}\Big)_{x_i},
 \label{par_cons}
 \end{equation}
where $\nu_{i,j}$ is a general diffusivity tensor. This linear problem is also to be seen as a model problem for nonlinear diffusion operators in which, for example, this tensor is derived by an algebraic closure of turbulence, see e.g. \cite{decoene:2009}, \cite{girard:1990}.
 It is therefore of interest to extend SL methods for diffusive problems to this kind of formulation, to avoid rewriting equation (\ref{par_cons}) as an advection--diffusion equation with a drift term that is due to the variable diffusion coefficient, rather than to the physical advection process.
 This would allow to address a number of accuracy and efficiency issues  concerning the discretization of advection--diffusion equation systems in
 a coherent SL framework. Indeed, in typical anisotropic meshes employed in environmental applications, the diffusion terms introduce a remarkable stiffness in the differential problem to be solved, which is usually handled by coupling an Eulerian discretization in space with an implicit time stepping method. This increases the computational cost per time step, requires more communication on parallel machines and can decrease accuracy due to the introduction of splitting errors.
With this direction of work in mind, we propose in this paper a general concept for the SL approximation of a linear second-order balance equation in divergence form based on a modified version of the scheme proposed in \cite{teixeira:1999}. We see this as an intermediate step towards the (ongoing) development of a fully conservative
SL method for advection-diffusion equations. This would be extremely appealing   since a number of mass conservative, flux-form SL methods have been proposed in the last couple of decades, see e.g. \cite{rancic:1995}, \cite{lin:1996}, \cite{leonard:1996}, \cite{frolkovic:2002}, \cite{nair:2002a}, \cite{zerroukat:2004}, \cite{restelli:2006}, \cite{qiu:2011a}, \cite{qiu:2011b}.  Although the proposed method is only first order in time, achieving higher order accuracy with unconditionally stable and robust schemes is not an easy task (see e.g. \cite{kalnay:1988}, \cite{teixeira:1999}, \cite{wood:2007}) and the approach we propose fits well in a coherent SL approach while achieving the same effective accuracy displayed by more standard discretizations.

  In section \ref{general}, the basic notation and definitions necessary to present SL methods will be introduced. In section \ref{second_ord}, our novel SL discretization is introduced, while an analysis of its consistency and stability properties will be presented in sections \ref{consistency} and \ref{stability}. Section \ref{ndim} treats the extension of the scheme to problems in general dimension and with nonlinear diffusivity. Finally, some numerical results obtained with the proposed method in both linear and nonlinear models will be presented in section \ref{tests}, while some conclusions on the potential advantages of our approach will be drawn in section \ref{conclu}.

\section{Basic facts and notations about semi-Lagrangian methods}
\label{general}
We collect in this section all the basic definitions and theoretical concepts needed in the paper. First, we briefly review the SL strategy to treat the advection equation
\begin{equation}\label{trasp}
\begin{cases}
      u_t+f(x,t) u_x=0,& (x,t)\in\R \times [0,T]\\
      u(x,0)=u_0(x)& x\in\R.\\
\end{cases}
\end{equation}
The construction of large time-step schemes for \eqref{trasp} stems from the application of the classical  method of characteristics.
 The system of characteristic curves $X(x,t;s)$ for \eqref{trasp} is defined by the solutions of:
\begin{equation}\label{caratt}
\begin{cases}
\displaystyle\frac{d}{ds}X(x,t;s)=f(X(x,t;s),s),\\
X(x,t;t)=x.
\end{cases}
\end{equation}
The solution of \eqref{trasp} is constant along such curves, which means that the following representation formula
\begin{equation}\label{rappres}
u(x,t)=u(X(x,t;t+\tau),t+\tau)
\end{equation}
holds for the solution $u$. Writing \eqref{rappres} with $\tau=-\Delta t$, we have the time-discrete version
\begin{equation}\label{upw}
u(x,t)=u(X(x,t;t-\Delta t),t-\Delta t).
\end{equation}

Discretizing the representation formula \eqref{rappres} we obtain the advective SL approximation. More precisely, we denote by $\Delta t$ and $\Delta x$ respectively the time and space discretization steps, with $t_n=n\Delta t$ for $n\in [0,T/\Delta t]$. The space grid is supposed to be infinite and uniform, that is, for $i\in\Z,$ $x_i= i\Delta x.$
We will denote the numerical solutions of \eqref{trasp} at time $t_n$ by the vector $V^n=(v_i^n)_{i\in\Z}$. Numerical solutions are usually analyzed in the $l^2$ norm, so that
\begin{equation*}
\|W\|_2 = \left(\Delta x \sum_{j\in\Z} w_j^2\right)^{1/2}.
\end{equation*}
The characteristics $X$ defined by \eqref{caratt} will be replaced by their numerical approximations $X^\Delta,$ whose specific form will not be discussed here. We will possibly use the shorthand notation
\begin{equation}\label{z_i}
z_i=X^\Delta(x_i,t^{n+1};t^n)
\end{equation}
to denote the foot of the approximate characteristic starting from $x_i$. For simplicity, this notation neglects the possible dependence of $z_i$ on $n$. We assume that the approximation $X^\Delta$ is consistent with order $r\ge 1$, that is,
\[
\left|X^\Delta(x,t;t-\Delta t)-X(x,t;t-\Delta t)\right| \le O\left(\Delta t^r\right).
\]
In advective SL schemes, \eqref{rappres} is discretized by replacing the exact upwinding $X$ with $X^\Delta$ and the value of $u$ at the foot of a characteristic with an interpolation $I_p$:
\begin{equation}\label{sl}
v_i^{n+1}=I_p[V^n] \left(X^\Delta(x_i,t^{n+1};t^n)\right)
\end{equation}
where $v_i^{n+1}$ is the approximation of $u(x_i,t^{n+1})$, and $I_p$ is an interpolation operator (e.g., a polynomial interpolation of degree $p$) which is assumed to satisfy the condition
\[
I_p[V](x_i) = v_i,
\]
and which, given a smooth function $w$ and a vector $W$ such that $w_i = w(x_i)$, satisfies
\[
\|w-I_p[W]\|_\infty \le O\left(\Delta x^{p+1}\right).
\]
We remark also that the consistency error for the schemes under consideration reads
\begin{equation}\label{cons_err}
L(\Delta x, \Delta t) = O\left(\Delta t^r + \frac{\Delta x^{p+1}}{\Delta t} \right).
\end{equation}

\section{Semi-Lagrangian methods for second order problems in flux form}
\label{second_ord}
We now introduce an approach to adapt advective form SL schemes to diffusion equations in divergence form. 
We consider first the the simpler case of the constant-coefficient advection--diffusion equation
\[
u_t + au_x = \nu u_{xx}.
\]
Give a time step of length $\Delta t,$ we introduce an advective displacement $\alpha=a\Delta t$
and a diffusive displacement defined, by standard scaling arguments, as
\[
\delta = \sqrt{2\Delta t\>\nu}.
\]
It is possible to construct an explicit approximation in the abstract form:
\begin{equation*}
u(x_i,t_{n+1}) \approx \frac{1}{2}u(x_i-\alpha+\delta,t_n) + \frac{1}{2}u(x_i-\alpha-\delta,t_n).
\end{equation*}
By Taylor expansion around the values $u(x_i-\alpha\pm\delta,t_n)$ it is easy to see  that
\begin{eqnarray*}
\frac{1}{2}u(x_i-\alpha+\delta,t_n) + \frac{1}{2}u(x_i-\alpha-\delta,t_n) & = & u(x_i,t_n) + \\
&& \hspace{-3cm} + \Delta t \left( -au_x + \nu u_{xx}(x_i,t_n) \right) + O(\Delta t^2),
\end{eqnarray*}
thus implying that a scheme based on this approach can only be of  first order with respect to $\Delta t$.
As far as the advective term is concerned, this approximation coincides with the standard advective form SL method,
while for the diffusive term this approximation strategy has a considerable literature and is traditionally derived from the Feynman--Kac stochastic representation formula for the solution of diffusion equations. We refer the reader to the papers cited in section \ref{intro} and to \cite{ferretti:2010} for a review of the related literature, as well as for a generalization to multidimensional problems and, partly, to higher order consistency rates.
 
To adapt this approach to more general problems, we start by examining the case of a variable-coefficient linear diffusion equation in divergence form:
\begin{equation}\label{diffusione}
u_t = (\nu(x,t) u_x)_x.
\end{equation}
The technique   proposed  here for treating equation \eqref{diffusione} is based on a simplified version of the scheme proposed in \cite{teixeira:1999}.
It consists in introducing, for each node, two possibly different displacements  $\delta_i^\pm$
to be suitably defined. This allows to define the SL update 
by the abstract   operator:
\[
u(x_i,t_{n+1}) \approx \frac{1}{2}u(x_i+\delta^+_i,t_n) + \frac{1}{2}u(x_i-\delta^-_i,t_n).
\]
In practice, a numerical interpolation is needed to recover the two values $u(x_i\pm\delta^\pm_i,t_n)$, so that the actual scheme takes the form
\begin{equation}\label{schema_tm}
v_i^{n+1} = \frac{1}{2}I[V^n](x_i+\delta^+_i) + \frac{1}{2}I[V^n](x_i-\delta^-_i).
\end{equation}
The displacements $\delta_i^\pm$ are given as the solution of the equation
\begin{equation}\label{delta_t}
\delta^\pm_i = \sqrt{2\Delta t\>\nu(x_i\pm\delta^\pm_i,t^n)}.
\end{equation}
It is to be remarked that this approach allows to increase easily the spatial
accuracy of the method by application of higher order interpolation operators, independently of the mesh
regularity. This contrasts with all standard finite difference and finite volume discretizations
of diffusion operators, for which second order accuracy is typically achieved only on uniform meshes.
Note also that \eqref{delta_t} is in fixed point form and, once denoted by $T^\pm_i$ the right-hand side of \eqref{delta_t}, it could be solved iteratively as
\begin{equation}\label{iteraz}
\delta^\pm_{i,k} = T^\pm_i\left(\delta^\pm_{i,k-1}\right),
\end{equation}
provided each $T^\pm_i$ is a contraction. On the other hand, an immediate computation shows that
\begin{equation}\label{Lip_T}
\left\|{T^\pm_i}'\right\|_\infty \le \sqrt{\frac{\Delta t}{2}}\> \frac{\|\nu_x\|_\infty}{\inf \nu},
\end{equation}
which implies that $T_i$ is contractive (at least for $\Delta t$ small enough) if $\nu$ is Lipschitz continuous and bounded from below by some positive constant. However, situations in which the diffusivity presents abrupt variations, or in which it becomes very small, could cause the iterative computation of \eqref{delta_t} to break down. In this case, bisection could be a less efficient but more robust choice.
It should also be observed that this approach to the computation of the diffusive displacement
is very similar to the standard fixed point technique proposed in \cite{robert:1981} and
widely employed in meteorological applications for the computation of advective displacements.

It is also to be remarked that,  if \eqref{iteraz} is used to solve \eqref{delta_t}, 
an a priori estimate of the accuracy of the computed solution is feasible.
Indeed, by a standard result on fixed-point iterations one has
\begin{equation}
\left| \delta^\pm_{i,k} - \delta^\pm_i \right| \le L_T^k \left|\delta^\pm_{i,0} - \delta^\pm_i \right|,
\end{equation}
where $L_T$ denotes the Lipschitz constant of $T^\pm_i$. Setting $L_T=C_L \Delta t^{1/2}$ (for some constant $C_L$) as provided by \eqref{Lip_T}, and taking as initial estimate
\[
\delta^\pm_{i,0} = \sqrt{2\Delta t\>\nu(x_i,t^n)}
\]
for which one has the initial error estimate
\[
\left|\delta^\pm_{i,0} - \delta^\pm_i \right| = O\left(\Delta t^{1/2}\right),
\]
we finally obtain
\begin{equation}\label{stima_iteraz}
\left| \delta^\pm_{i,k} - \delta^\pm_i \right| = O\left(\Delta t^\frac{k+1}{2}\right).
\end{equation}
We will use this estimate later on in the consistency analysis of the scheme.

The previous approach can be extended to the discretization of the  linear second-order equation
\begin{equation}\label{diftra_adv}
\begin{cases}
u_t + f(x,t)u_x=(\nu(x,t) u_x)_x & (x,t)\in\R\times [0,T]\\
u(x,0) = u_0(x) & x\in\R
\end{cases}
\end{equation}
in which both the advection and the diffusion terms appear. In this case, advective displacements 
$\alpha_i = x_i-z_i$ (with $z_i$ defined by \eqref{z_i}) and diffusive displacements $\delta^\pm_{i} $ have to be computed. The former displacement can be computed by standard fixed point iterations as proposed in \cite{robert:1981} or by the sub-stepping approaches described, e.g., in \cite{giraldo:1999}, \cite{rosatti:2005}.
In the present approach, these two steps are performed independently and the total displacement is simply
obtained by adding the advective and the diffusive displacement. This separate computation entails a first order splitting error, which is compatible with the overall accuracy of the proposed method. Higher order displacement computation procedures could be devised, but they would cause an increase of the method's complexity and computational cost. The resulting advective SL method for equation \eqref{diftra_adv} can thus be written as
\begin{equation}\label{schema_tm_adv}
v_i^{n+1} = \frac{1}{2}I[V^n](z_i+\delta^+_i) + \frac{1}{2}I[V^n](z_i-\delta^-_i).
\end{equation}
Note that for a vanishing viscosity $\nu$ the scheme reduces, without any loss in stability properties, to the corresponding scheme \eqref{sl} for the inviscid problem. Note also that, for large values of the time step, the two points $z_i\pm\delta^\pm_i$ could be relatively far from each another. This situation has no consequence when handling smooth solutions, but may cause the smaller scales to be severely underresolved, at least in the first steps, when nonsmooth solutions are considered, as analysed in \cite{ferretti:2010}.  Adaptive time-stepping strategies to reduce this drawback have been studied (e.g., in \cite{carlini2006time} for the case of the Mean Curvature equation).

\section{Consistency}
\label{consistency}
For the consistency analysis, we will restrict to the approximation of the pure diffusion equation \eqref{diffusione}.
For the advection terms, the reader may refer to
the results reported in \cite{falcone:2013}. We start by rewriting the method in abstract form as a linear (in fact, convex) combination of pointwise values:
\begin{equation}\label{schema}
u(x_i,t_{n+1}) \approx A_i^+ u(x_i+\delta_i^+,t_n) + A_i^- u(x_i-\delta_i^-,t_n).
\end{equation}
Expressing in \eqref{schema} the values $u(x_i\pm\delta_i^\pm,t_n)$ by means of their Taylor expansions centered at $x_i$, we get
\begin{equation}
u(x_i\pm\delta_i^\pm,t_n) = u \pm \delta_i^\pm u_x + \frac{{\delta_i^\pm}^2}{2} u_{xx} \pm \frac{{\delta_i^\pm}^3}{3!} u_{xxx} + O\left({\delta_i^\pm}^4\right)
\end{equation}
(where $u,\ldots,u_{xxx}$ are computed at $(x_i,t_n)$), which yield, when used in \eqref{schema},
\begin{eqnarray}\label{sviluppo}
&&\hspace{-2cm}A_i^+ u(x_i+\delta_i^+,t_n) + A_i^- u(x_i-\delta_i^-,t_n) = \nonumber \\
& = & \left(A_i^++A_i^-\right)u + \left(A_i^+\delta_i^+-A_i^-\delta_i^-\right) u_x + \nonumber \\
&& +\frac{1}{2} \left(A_i^+{\delta_i^+}^2+A_i^-{\delta_i^-}^2\right) u_{xx} + \nonumber \\
&& +\frac{1}{3!} \left(A_i^+{\delta_i^+}^3-A_i^-{\delta_i^-}^3\right) u_{xxx} + O\left({\delta_i^\pm}^4\right).
\end{eqnarray}
On the other hand, for  equation \eqref{diffusione} we have:
$$
u_t =  (\nu u_x)_x  =  \nu_xu_x + \nu u_{xx},
$$
that is, expanding to first order in time,
\[
u(x_i,t_{n+1}) = u + \Delta t (\nu_x u_x + \nu u_{xx}) + O\left(\Delta t^2\right).
\]
Comparing this expression with \eqref{sviluppo}, we obtain a (first-order) consistent approximation under the conditions
\begin{equation}\label{sistema}
\begin{cases}
A_i^+ + A_i^- = 1 + O(\Delta t^2) \\
A_i^+\delta_i^+ - A_i^-\delta_i^- = \Delta t \> \nu_x + O(\Delta t^2) \\
A_i^+{\delta_i^+}^2 + A_i^-{\delta_i^-}^2 = 2\Delta t \> \nu + O(\Delta t^2) \\
A_i^+{\delta_i^+}^3 - A_i^-{\delta_i^-}^3 = O(\Delta t^2).
\end{cases}
\end{equation}
Note that, in fact, the additional condition ${\delta_i^\pm}^4 = O\left(\Delta t^2\right)$ should also be added, but, according 
to  definition \eqref{delta_t}, this is satisfied automatically.
We apply now the abstract conditions \eqref{sistema} to prove consistency of the proposed scheme.

In the advective form SL method \eqref{schema_tm}, we have $A_i^+=A_i^-=1/2$, so that the first condition in \eqref{sistema} is satisfied. The displacements $\delta_i^\pm$ are defined by \eqref{delta_t} from which, using a Taylor expansion for $\nu$ and squaring both sides, we get
\begin{equation}\label{quadrato}
{\delta_i^\pm}^2 = 2\Delta t \> \left(\nu\pm\nu_x \delta_i^\pm + \frac{1}{2}\nu_{xx}{\delta_i^\pm}^2+ O(\Delta t^{3/2})\right).
\end{equation}
Concerning now the second condition in \eqref{sistema}, we can rewrite its left-hand side via \eqref{quadrato} as
\begin{eqnarray}\label{seconda}
\frac{1}{2}(\delta_i^+ - \delta_i^-) & = & \frac{{\delta_i^+}^2 - {\delta_i^-}^2}{2\left(\delta_i^+ + \delta_i^- \right)} = \nonumber \\
& = & \frac{2\Delta t \left(\delta_i^+ + \delta_i^- \right)\nu_x + O(\Delta t^{5/2})}{2\left(\delta_i^+ + \delta_i^- \right)} = \nonumber \\
& = & \Delta t \> \nu_x + O(\Delta t^2).
\end{eqnarray}
Using again \eqref{quadrato} along with \eqref{seconda}, we also obtain:
\begin{equation*}
{\delta_i^+}^2 + {\delta_i^-}^2 = 4\nu\Delta t +O(\Delta t^2)
\end{equation*}
so that the third condition is also satisfied. Finally, the left-hand side in the fourth condition may be written as
\begin{equation*}
\frac{1}{2}\left({\delta_i^+}^3 - {\delta_i^-}^3\right) = \frac{1}{2}\left(\delta_i^+ - \delta_i^- \right)\left({\delta_i^+}^2 + \delta_i^+\delta_i^- + {\delta_i^-}^2\right),
\end{equation*}
where, using again \eqref{seconda}, it is immediate to see that the right-hand side is $O(\Delta t^2)$, and the fourth condition is satisfied.
It can be observed that, thanks to  \eqref{stima_iteraz}, if the displacement are computed by fixed point iterations,
three such iterations suffice to achieve an $O(\Delta t^2)$ error, which is enough to preserve the consistency order of the scheme.
Finally, introducing a space discretization as outlined in Sec. \ref{general}, we note that the
 pointwise values $u(x_i\pm\delta_i^\pm)$ are approximated with accuracy of order $O(\Delta x^{p+1})$. The estimate \eqref{cons_err} then holds with $r=1$.

\section{Stability}
\label{stability}

To discuss the stability of the proposed scheme, we  show that it can be recast as a convex combination of advective SL
schemes. In this framework, we can easily prove stability for constant coefficient equations, while variable coefficient equations require a deeper study.

The advective SL method can be rewritten in matrix form as
\begin{equation}\label{media}
V^{n+1} = \frac{1}{2} B^+ V^n + \frac{1}{2} B^- V^n,
\end{equation}
where, for example, the matrix $B^+$ represents the interpolation step
\begin{equation}\label{inviscid_tm}
w_i^{n+1} = I[W^n](x_i+\delta_i^+)
\end{equation}
and $B^-$ plays the same role for a scheme with upwinding $-\delta_i^-$. Now, in the constant-coefficient case,
 we have $\delta_i^\pm \equiv \delta$, and it is known that \eqref{inviscid_tm} is nonexpansive in the 2-norm 
 for any order of Lagrange interpolation, that is,
\[
\left\| B^\pm \right\|_2 = 1.
\]
Therefore, the complete scheme \eqref{schema_tm} is also stable in the same norm.
Notice that the technique used in \cite{ferretti:2013a} to prove stability of SL schemes for the variable coefficient case  seems unsuitable to treat \eqref{inviscid_tm}. Indeed, since the displacements $\delta_i^\pm$ are of the order of $\Delta t^{1/2}$, a straightforward application of the ideas of \cite{ferretti:2013a} to the present cas  would provide a perturbation of order $O(\Delta t^{1/2})$  of the equivalent Lagrange-Galerkin scheme, which does not guarantee stability.

On the other hand, for at least two situations a complete stability analysis is feasible. The first occurs if the advecting field $f$ depends on $x$, but $\nu$ does not. In this case, the $\delta_i^\pm$ only introduce a constant displacement and the stability analysis of \cite{ferretti:2013a,ferretti:2013b} remains unchanged. The second occurs in case both $f$ and $\nu$ depend on $x$, but the scheme \eqref{inviscid_tm} is monotone (e.g., with $\mathbb{P}_1$ interpolation). In this case, the resulting scheme for the diffusion equation is also monotone, as it is easy to check. 

 \section{Generalization to multiple space dimensions and nonlinear problems}
\label{ndim}
In this section, we briefly discuss the extension of the proposed approach to the $d$-dimensional case. 
The extension is relatively straightforward if the diffusivity matrix has the diagonal form
\begin{equation}\label{diff_diag}
\Lambda(x,t) = \text{diag} (\nu_1(x,t),\ldots,\nu_d(x,t)),
\end{equation}
and in particular, for a variable but isotropic diffusion (for which $\nu_1(x,t)=\cdots =\nu_d(x,t)$). Then, the diffusion equation reads
\begin{eqnarray*}
u_t & = & \text{div}(\Lambda(x,t)\nabla u) = \\
& = & \sum_{j=1}^d (\nu_j(x,t)u_{x_j})_{x_j}.
\end{eqnarray*}
Performing a first-order expansion of the solution with respect to time, and rearranging the various terms, we get
\begin{eqnarray*}
u(x,t+\Delta t) & = & u(x,t) + \Delta t \sum_{j=1}^d (\nu_j u_{x_j})_{x_j}(x,t) + O(\Delta t^2) = \\
& = & \sum_{j=1}^d \left[ \frac{u(x,t)}{d}+\frac{\Delta t}{d} (d\nu_j u_{x_j})_{x_j}(x,t)\right] + O(\Delta t^2) = \\
& = & \frac{1}{d} \sum_{j=1}^d \left[ u(x,t) + \Delta t (d\nu_j u_{x_j})_{x_j}(x,t)\right] + O(\Delta t^2),
\end{eqnarray*}
which shows that, up to first-order accuracy, the $d$-dimensional version can be obtained by averaging the diffusion operators in each direction, with the only modification of scaling each one-dimensional diffusivity by a factor $d$. For example, the 2-dimensional equation
\[
u_t = (\nu_1(x_1,x_2,t)u_{x_1})_{x_1} + (\nu_2(x_1,x_2,t)u_{x_2})_{x_2}
\]
would be approximated by the scheme (written with an obvious notation at the node $x_i=(x_{1,i},x_{2,i})$)
\begin{eqnarray*}
v_i^{n+1} & = & \frac{1}{4} \left( I[V^n](x_{1,i}+\delta^+_{1,i},x_{2,i}) + I[V^n](x_{1,i}-\delta^-_{1,i},x_{2,i}) + \right.\\
&& \left. + I[V^n](x_{1,i},x_{2,i}+\delta^+_{2,i}) + I[V^n](x_{1,i},x_{2,i}-\delta^-_{2,i}) \right),
\end{eqnarray*}
in which the displacements $\delta^\pm_{j,i}$ are defined by
\begin{eqnarray*}
&& \delta^\pm_{1,i} = \sqrt{4\Delta t\>\nu_1(x_{1,i}\pm\delta^\pm_{1,i},x_{2,i},t^n)} \\
&& \delta^\pm_{2,i} = \sqrt{4\Delta t\>\nu_2(x_{1,i},x_{2,i}\pm\delta^\pm_{2,i},t^n)}.
\end{eqnarray*}
To treat the more general case
\[
u_t = \text{div}(A(x,t)\nabla u)
\]
for a symmetric, positive semidefinite matrix $A(x,t)$, we exploit the existence of transformed coordinates $\xi, $ 
that for given $\bar x$ and $\bar t$ we will denote by $\xi=\Theta(\bar x,\bar t) x$ (with $\Theta(\bar x,\bar t)\in\R^{d\times d}$), such that
\[
\Lambda(\xi, t)=\Theta(\bar x,\bar t) A(x, t) \Theta(\bar x,\bar t)^{-1}
\]
is diagonal at $\bar\xi = \Theta(\bar x,\bar t) \bar x$. In this local system of coordinates, the diffusion operator is rewritten at $(\bar x,\bar t)$ as
\begin{eqnarray*}
\text{div}_x(A(x, t)\nabla_x u) & = & \text{div}_\xi\left(\Theta(\bar x,\bar t) A(x, t)\Theta(\bar x,\bar t)^{-1} \nabla_\xi u\right) = \\
& = & \sum_{j=1}^d \left(\nu_j(\xi, t\right)u_{\xi_j})_{\xi_j},
\end{eqnarray*}
and the problem is brought back to the diagonal case by working in the variables $\xi$. Note, however, that the transformation $\Theta$ should be recomputed at each point $(x_i,t_n)$.

Finally, we add some remark on the application of the scheme to nonlinear problems in the form
\[
u_t = \text{div}(A(x,t,u,\nabla u)\nabla u).
\]
In all the analysis above, we have assumed that the dependence on $t$ of the diffusivity tensor $A$ is treated by freezing the diffusivity at the previous time step. 
This amounts to treat the nonlinear problem, when advancing from $t_n$ to $t_{n+1}$, by a linearization in the form
\[
u_t = \text{div}(A(x,t_n,u(t_n),\nabla u(t_n))\nabla u)
\]
(note that, whenever diffusivity depends on the solution, it is only known at the nodes and needs to be interpolated to compute the displacements $\delta^\pm_{j,i}$).
This extension clearly requires some additional caution. We will present a couple of numerical examples in the next section, trying to point out some specific devices to make the approach work correctly.

 \section{Numerical experiments}
\label{tests}
Several numerical experiments have been carried out with simple implementations of the advective SL method proposed above,
in order to assess its accuracy and stability features also in more complex cases than those allowing a complete theoretical analysis. The accuracy of the proposed discretization has been evaluated against reference solutions obtained by alternative discretizations in space and time.

\subsection{Constant coefficient case}
\label{constcoeff}
In a first set of numerical experiments, the constant coefficient
advection equation
\[
u_t + au_x = \nu u_{xx} \ \ \ \ \ x\in[0,L]
\]
was considered, on an interval $[0,L] $ with $L=10.$  Periodic boundary conditions
were assumed and a Gaussian profile centered at $L/2$ was considered as the initial condition.
In this case, the exact solution  can be computed up to machine accuracy by separation
of variables and computation of its Fourier coefficients on a discrete mesh of
$N$ points with spacing $\Delta x=L/N.$ We consider the SL methods
described in the previous sections  on a time interval $[0,T] $ with $T=2.75 $
with time steps defined as $\Delta t=T/M.$ The Courant number and the stability
parameter of standard explicit discretization of the diffusion operator are defined
as $C=a\Delta t / \Delta x $ and $\mu= \nu \Delta t/ 2\Delta x^2, $ respectively. We consider the case
with $a=0$, $\nu=0.05 $ first, whose results are reported in tables \ref{table1}-\ref{table2}, 
for the SL method with linear and cubic reconstructions, respectively.
The parallel results for the case with $a=1, \nu=0.05$ are reported in tables \ref{table3}-\ref{table4}.

It can be observed that, as expected, the SL method has much greater accuracy with cubic interpolation
than with linear interpolation, allowing to achieve almost second order convergence
rates and to improve the accuracy in the limit of larger time step sizes, as typical with SL schemes.
In order to carry out a  comparison with a more standard technique, the same test has also
been run employing a fourth order centered finite difference discretization
for advection, conservative second order finite differences for the diffusion term and an off-centered
Crank Nicolson method (i.e., the so called $\theta-$method, with $\theta=0.52$). This combination, that
is quite widely employed in many environmental applications especially for the diffusion terms,
yields a spatial accuracy for the advection term that should be comparable with that
of SL with cubic interpolation. Also the time discretization is first order, as the SL method
proposed here for the diffusion terms. From the results of this test, reported in table
\ref{table5}, it can be seen that the proposed SL technique is superior to the standard
discretization approach at all resolutions. Equivalent results have instead
been obtained with  the two approaches  
 if $\theta=0.5 $ was chosen, which is however never done in realistic 
applications. Finally, the SL method with cubic interpolation and the just described
finite difference discretization have been  compared also for larger values
of the stability parameters  $C$ and $\mu.$ Results reported in table \ref{table6} show that
the SL method is clearly superior in accuracy also in the limit of large time steps.

 \begin{table}[t]
\centering
\begin{tabular}{||c|c|c|c|c||} 
\hline
\multicolumn{3}{||c|}{Resolution} & \multicolumn{2}{c||}{Relative error} \\
\hline
$N$ & $M$ &$\mu $ & $l_2$ &  $l_{\infty}$ \\
\hline \hline
$200$ &  $100$ & $0.55$  & $7.70\cdot 10^{-3}$ & $6.60\cdot 10^{-3}$  \\
\hline
$200$ & $ 200 $ & $0.138$ & $5.68\cdot 10^{-2}$ & $4.99\cdot 10^{-2}$ \\
\hline
$400$ & $ 100 $ & $1.1$ & $3.80\cdot 10^{-3}$ & $3.30\cdot 10^{-3}$  \\
\hline

$400$ & $ 200 $ & $0.275$ & $1.71\cdot 10^{-2}$ & $1.97\cdot 10^{-2}$   \\
\hline
\end{tabular}
\caption{Errors for constant coefficient diffusion equation, SL method with linear interpolation.}
\label{table1} 
 \end{table}



 \begin{table}[t]
 \centering
\begin{tabular}{||c|c|c|c|c||} 
\hline
\multicolumn{3}{||c|}{Resolution} & \multicolumn{2}{c||}{Relative error} \\
\hline
$N$ & $M$ &$\mu $ & $l_2$ &  $l_{\infty}$ \\
\hline \hline
$200$ &  $100$ & $0.55$ & $2.39\cdot 10^{-3}$ & $2.43\cdot 10^{-3}$    \\
\hline
$200$ & $200 $ & $0.138$ & $4.75\cdot 10^{-4}$ & $1.20\cdot 10^{-4}$   \\
\hline
$400$ & $100 $ &  $1.1$ & $2.62\cdot 10^{-4}$ & $2.92\cdot 10^{-4}$    \\
\hline
$400$ & $200$ &  $0.275$ & $1.48\cdot 10^{-4}$ & $2.12\cdot 10^{-4}$     \\
\hline
\end{tabular}
\caption{Errors for constant coefficient diffusion equation, SL method with cubic interpolation.}
\label{table2} 
 \end{table}


  \begin{table}[t]
  \centering
\begin{tabular}{||c|c|c|c|c|c||} 
\hline
\multicolumn{4}{||c|}{Resolution} & \multicolumn{2}{c||}{Relative error} \\
\hline
$N$ & $M$ & $C$ & $\mu $ & $l_2$ &  $l_{\infty}$ \\
\hline \hline
$200$ & $100 $ & $0.55$ & $0.275$ & $3.28\cdot 10^{-2}$ & $3.76\cdot 10^{-3}$    \\
\hline
$200$ & $200 $ & $0.275$ & $0.138$ & $3.52\cdot 10^{-2}$ & $4.04\cdot 10^{-2}$   \\
\hline
$400$ & $100 $ & $1.1$ & $1.1$ & $3.00\cdot 10^{-3}$ & $3.06\cdot 10^{-3}$    \\
\hline
$400$ & $200 $ & $0.55$ & $0.55$ & $3.40\cdot 10^{-3}$ & $3.39\cdot 10^{-3}$     \\
\hline
\end{tabular}
\caption{Errors for constant coefficient advection--diffusion equation, SL method with linear interpolation.}
\label{table3} 
 \end{table}


 \begin{table}[t]
 \centering
\begin{tabular}{||c|c|c|c|c|c||} 
\hline
\multicolumn{4}{||c|}{Resolution} & \multicolumn{2}{c||}{Relative error} \\
\hline
$N$ & $M$ & $C$ & $\mu $ & $l_2$ &  $l_{\infty}$ \\
\hline \hline
$200$ & $100 $ & $0.55$ & $0.275$ & $4.92\cdot 10^{-4}$ & $9.70\cdot 10^{-4}$    \\
\hline
$200$ & $200 $ & $0.275$ & $0.138$ & $2.94\cdot 10^{-4}$ & $6.02\cdot 10^{-4}$   \\
\hline
$400$ & $100 $ & $1.1$ &  $1.1$ & $2.72\cdot 10^{-4}$ & $3.41\cdot 10^{-4}$    \\
\hline
$400$ & $200 $ & $0.55$ &  $0.55$ & $1.78\cdot 10^{-4}$ & $2.34\cdot 10^{-4}$     \\
\hline
\end{tabular}
\caption{Errors for constant coefficient advection--diffusion equation, SL method with cubic interpolation.}
 \label{table4}   
 \end{table}

 

  \begin{table}[t]
 \centering
\begin{tabular}{||c|c|c|c|c|c||} 
\hline
\multicolumn{4}{||c|}{Resolution} & \multicolumn{2}{c||}{Relative error} \\
\hline
$N$ & $M$ & $C$ & $\mu $ & $l_2$ &  $l_{\infty}$ \\
\hline \hline
$200$ & $100 $ & $0.55$ & $0.275$ & $1.6\cdot 10^{-3}$ & $1.90\cdot 10^{-3}$    \\
\hline
$200$ & $200 $ & $0.275$ & $0.138$ & $8.55\cdot 10^{-4}$ & $7.55\cdot 10^{-4}$   \\
\hline
$400$ & $100 $ & $1.1$ &  $1.1$ & $1.7\cdot 10^{-3}$ & $1.90\cdot 10^{-3}$  \\
\hline
$400$ & $200 $ & $0.55$ &  $0.55$ & $8.23\cdot 10^{-4}$ & $9.46\cdot 10^{-4}$    \\
\hline
\end{tabular}
\caption{Errors for constant coefficient advection--diffusion equation, off-centered Crank--Nicolson time discretization and standard finite difference spatial discretization.}
 \label{table5}   
 \end{table}


  \begin{table}[t]
 \centering
\begin{tabular}{||l|c|c|c|c||} 
\hline 
\multicolumn{3}{||c|}{Discretization} & \multicolumn{2}{c||}{Relative error} \\
\hline
 Scheme & $C $ & $\mu $ &$l_2$ &  $l_{\infty}$  \\
\hline \hline
 SL & 1.375 & 0.687 & $6.69\cdot 10^{-4}$ & $8.74\cdot 10^{-4}$   \\
 \hline 
FD+CN & 1.375 & 0.687 & $4.50\cdot 10^{-3}$ & $5.5\cdot 10^{-3}$   \\
\hline
SL & 2.75 & 1.375 & $1.4\cdot 10^{-3}$ & $1.7\cdot 10^{-3}$   \\
 \hline 
FD+CN & 2.75& 1.375 & $1.15\cdot 10^{-2}$ & $1.35\cdot 10^{-2}$   \\
\hline
SL & 5.5 & 2.75 & $2.8\cdot 10^{-3}$ & $3.4\cdot 10^{-3}$   \\
 \hline 
FD+CN & 5.5 & 2.75 & $3.60\cdot 10^{-2}$ & $4.04\cdot 10^{-2}$   \\
\hline

\end{tabular}
\caption{Errors for constant coefficient advection--diffusion equation, comparison of full SL discretization
versus off-centered Crank--Nicolson time discretization and standard finite difference
spatial discretization for large time step values.}
 \label{table6}   
 \end{table}

%

\begin{figure}
\begin{center}
\includegraphics[width=11cm]{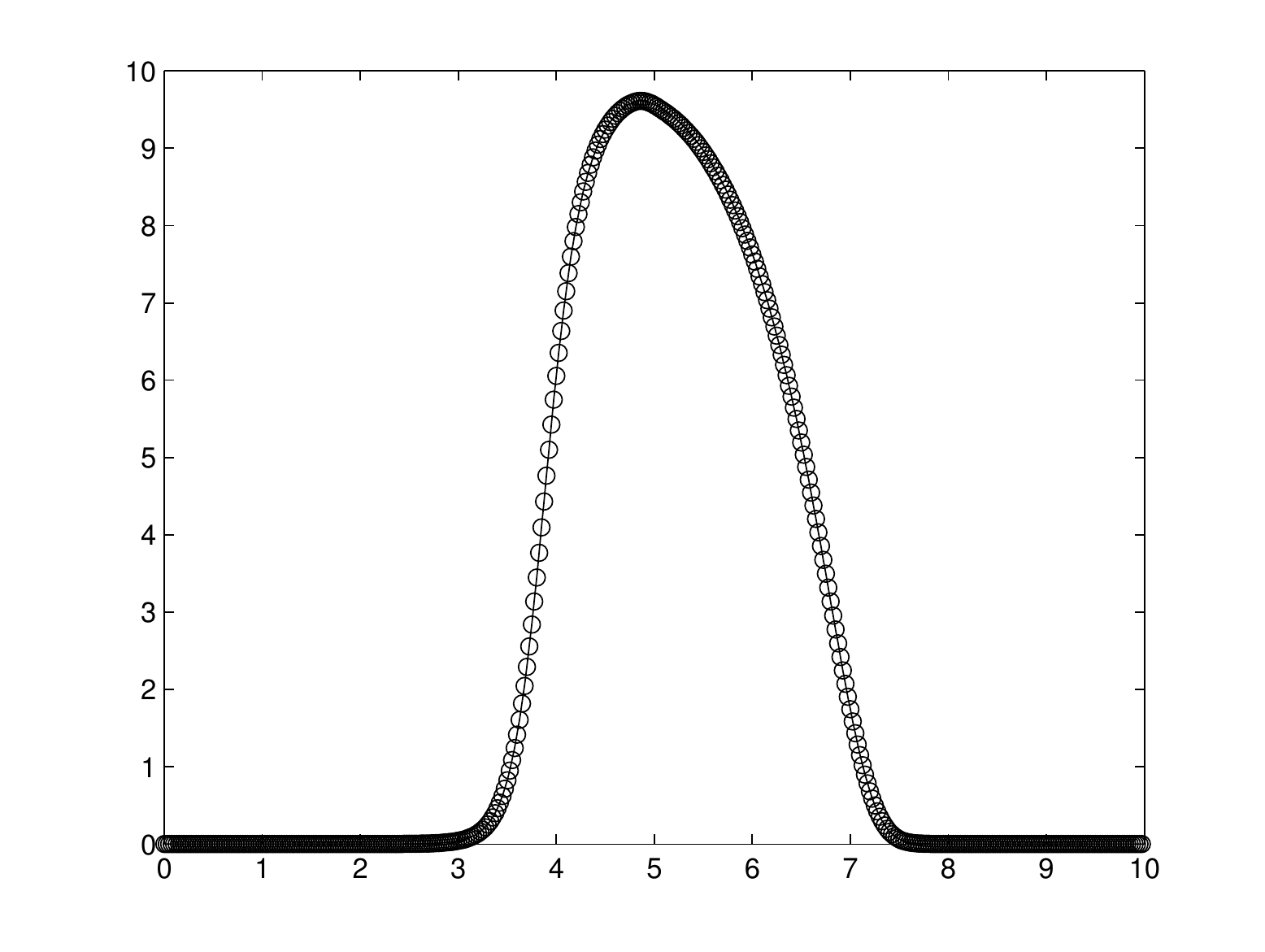}
\caption{One-dimensional case with variable coefficients: numerical solution by SL method (circles), reference solution (continuous line)}
\label{fig:vcoef1d}
\end{center}
\end{figure}

\subsection{Variable coefficient case: one space dimension}
\label{varcoeff1}

\[
u_t + f(x,t)u_x = (\nu(x,t) u_{x} )_x\ \ \ \ \ x\in[0,L]
\]
was then considered, on an interval $[0,L] $ with $L=10$  and on the 
 time interval $[0,T] $ with $T=4$
with time steps defined as $\Delta t=T/M.$ 
 Periodic boundary conditions
were assumed and gaussian profile centered at $L/3$ was considered as the initial condition.
The velocity and diffusivity field were given by
$$
f(x,t)= \frac 12+ \frac 12 \cos\left(\frac{6\pi x}L\right)\cos\left(\frac{2\pi t}T\right)
$$
$$
\nu(x,t)= \frac 1{100}+ \frac 1{25} \xi(x)\sin\left(\frac{2\pi t}T\right)^2 ,
$$
respectively, where $ \xi(x) $ denotes the characteristic function of the interval $[0.5L,0.8L].$ This choice highlights the possibility to use the
proposed method seamlessly also with strongly varying diffusion coefficients.
In this case, no exact solution is available and reference solutions were computed using
the finite difference method described in the previous section with a four times higher spatial
resolution, coupled to a high order  multistep stiff solver
in which a small tolerance and maximum time step value were enforced. 
A plot of the solution and reference solution at the
final time $T$ is displayed in figure \ref{fig:vcoef1d}.

A more quantitative assessment of the SL solution accuracy can be gathered from table \ref{vcoef_tab}.
It can be observed that, due to the time dependence of the diffusion coefficient and the
intrinsically first order nature of the time discretization proposed, the errors are higher than in the
corresponding constant coefficient case. On the other hand, they are reasonably small
with respect to the needs of many environmental applications and quite insensitive to the
choice of the time step.

  \begin{table}[t]
 \centering
\begin{tabular}{||c|c|c|c|c|c||} 
\hline
\multicolumn{4}{||c|}{Resolution} & \multicolumn{2}{c||}{Relative error} \\
\hline
$N$ & $M$ & $C$ & $\mu $ & $l_2$ &  $l_{\infty}$ \\
\hline \hline
$200$ & $10 $ & 4 & 0.8 & $1.2\cdot 10^{-2}$ & $1.92\cdot 10^{-2}$   \\
\hline
$200$ & $20 $ & 4 & 0.4 & $6.1\cdot 10^{-3}$ & $9.8\cdot 10^{-3}$   \\
\hline
$200$ & $100 $ & 0.8 & 0.08 & $1.3\cdot 10^{-3}$ & $2.30\cdot 10^{-3}$    \\
\hline
$200$ & $200 $ & 0.4 & 0.04 & $3.2\cdot 10^{-3}$ & $3.6\cdot 10^{-3}$   \\
\hline
$400$ & $100 $ & 1.6 & 0.32 & $1.4\cdot 10^{-3}$ & $2.0\cdot 10^{-3}$    \\
\hline
$400$ & $200 $ & 0.8 &  0.16 & $1.2\cdot 10^{-3}$ & $2.2\cdot 10^{-3}$     \\
\hline
\end{tabular}
\caption{Errors for variable coefficient advection--diffusion equation.}
 \label{vcoef_tab}   
 \end{table}


\subsection{Variable coefficient case: two space dimensions}
\label{varcoeff2}

\begin{figure}
\begin{center}

\includegraphics[width=11cm]{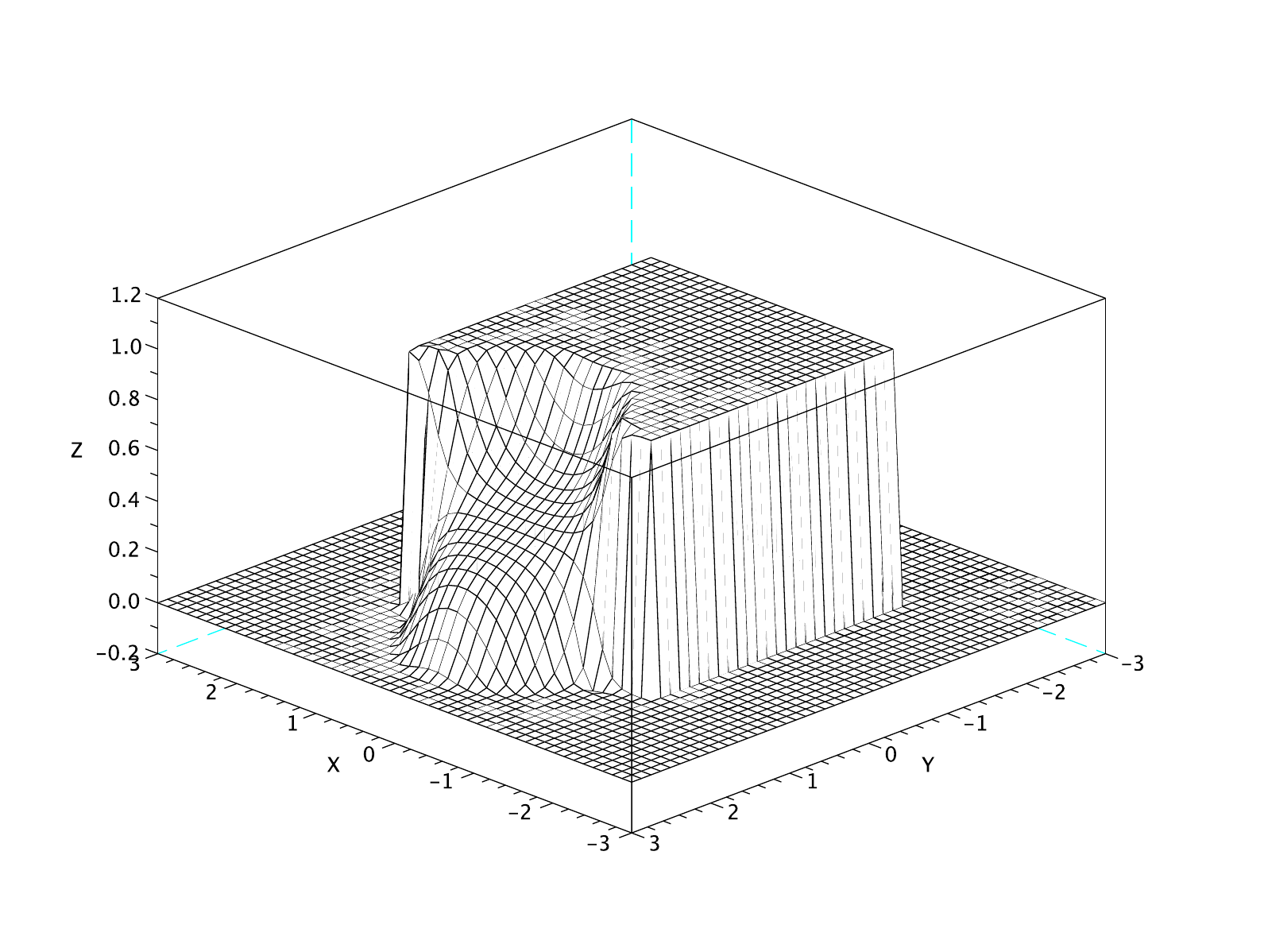}

\caption{Two-dimensional case with isotropic diffusivity}
\label{fig:isot2}
\end{center}
\end{figure}

In a first two dimensional  test, we consider the equation
\[
u_t = \text{div}(\nu(x)\nabla u)
\]
on $\Omega=[-3,3]^2$ with periodic boundary conditions. 
The initial condition given by the characteristic function of the set $\Sigma=[-1.5,1.5]^2$. The isotropic diffusivity $\nu$ is given by
\[
\nu(x) = e^{-5|x-x_0|^2},
\]
where $x_0=(1.5,0)$; the diffusion is therefore concentrated in a small region on the boundary of the set $\Sigma$. The effect of this diffusion is to move mass from the interior of $\Sigma$ to the exterior, in the neighbourhood of the point $x_0$. Fig. \ref{fig:isot2} shows the numerical solution at $T=1$, with a $50\times 50$ space grid, cubic interpolation and time step $\Delta t =0.05$.

\begin{figure}
\begin{center}
\includegraphics[width=0.65\textwidth]{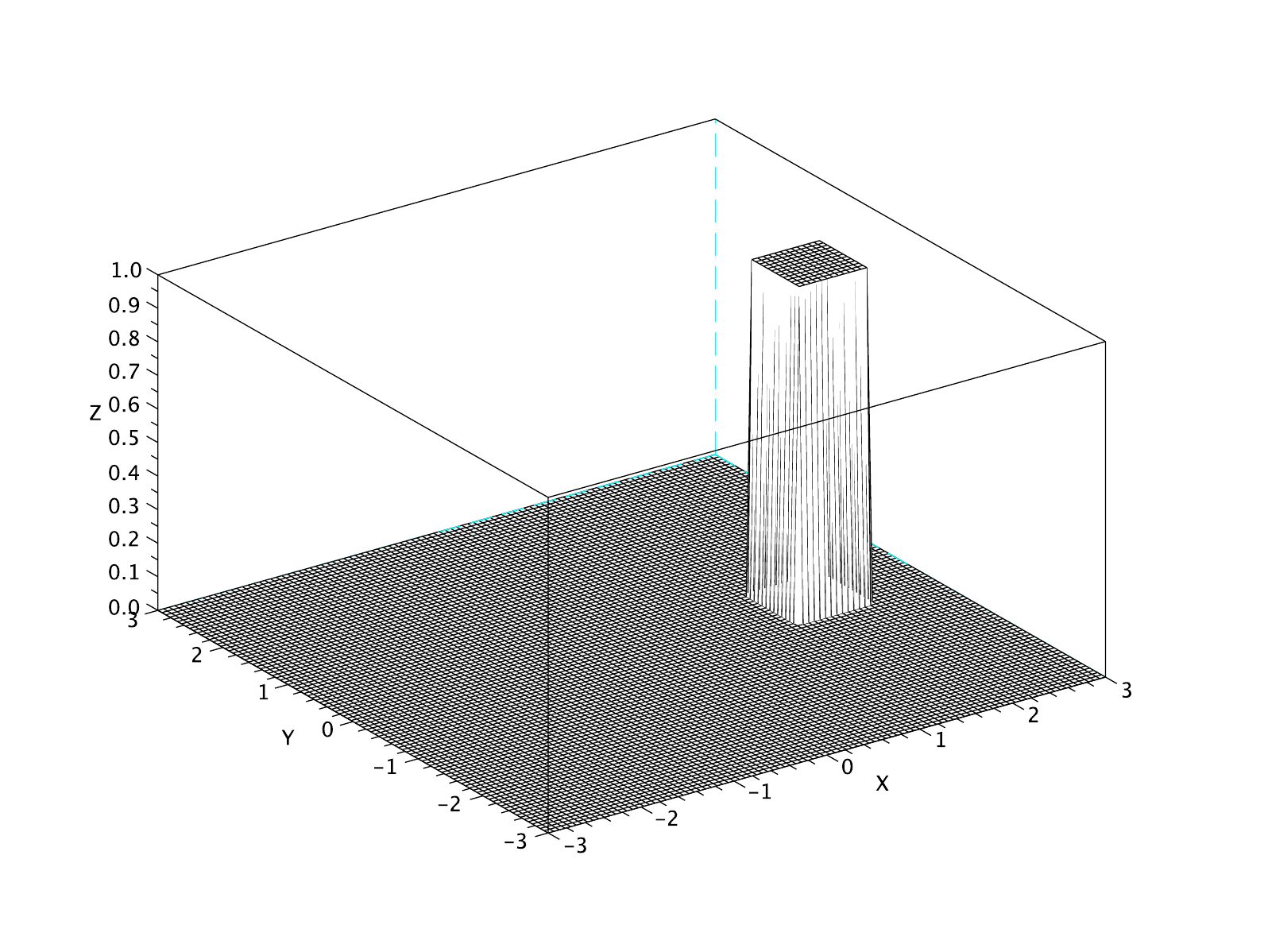}
\includegraphics[width=0.65\textwidth]{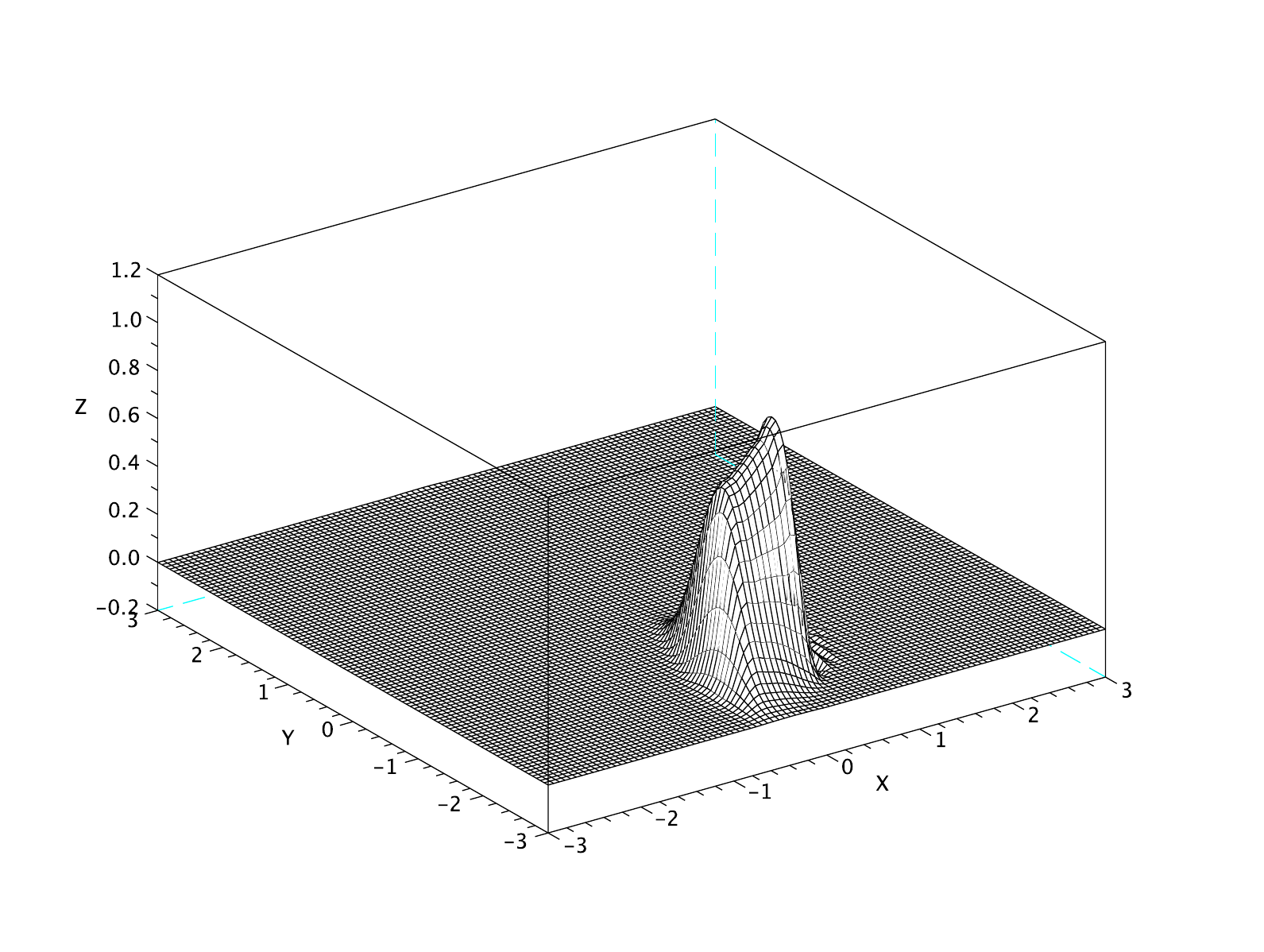}
\includegraphics[width=0.65\textwidth]{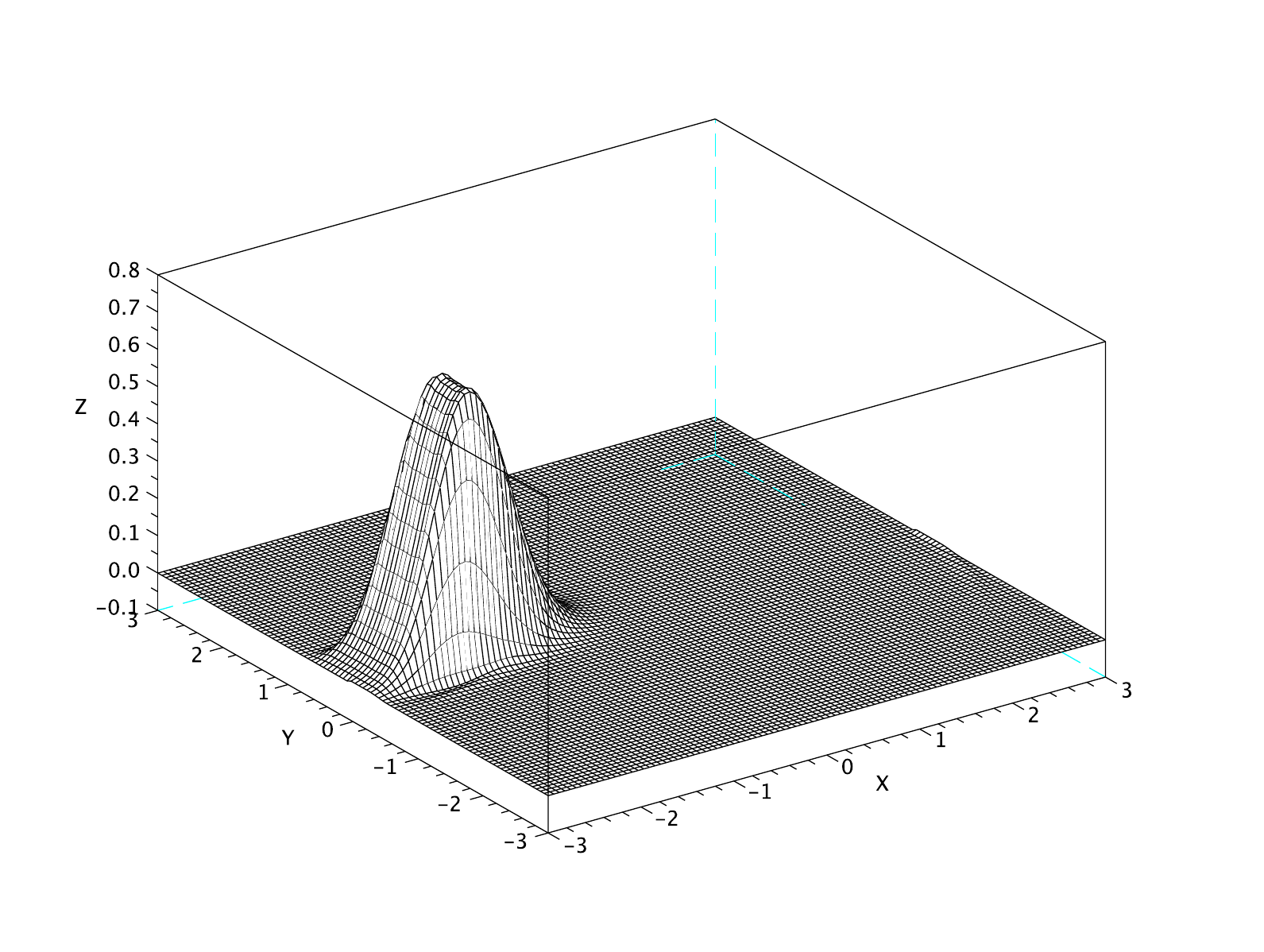}
\caption{Two-dimensional case with rotation and anisotropic diffusivity}
\label{fig:anisot2}
\end{center}
\end{figure}

In the second test, we consider the advection--diffusion equation
\[
u_t +f(x)\cdot \nabla u= \text{div}(\nu(x)\nabla u)
\]
with $f(x) =(x_2, -x_1 )$ and with spatially variable and anisotropic diffusivity of the form \eqref{diff_diag}, with
\[
\Lambda(x) = \begin{pmatrix}
0 & 0 \\
0 & e^{-5 x_1^2} \\
\end{pmatrix}.
\]
on the same domain $\Omega $ as in the previous case.
The initial condition is given by the characteristic function of a decentered square. The initial scalar field is then advected on circular trajectories, and diffused along the $x_2$ direction when crossing the axis $x_1=0$ (note that high values of the diffusivity are concentrated in a neighbourhood of this axis).
 Figure \ref{fig:anisot2} shows three snapshots of the solution, computed for $t\in[0,0.75]$ on a $100\times 100$ mesh, usingcubic interpolation and $\Delta t=0.00625$. Note that the plots have different scales for the $x_3$-axis.

\subsection{Gas flow in porous media}
\label{porous}
As  fist nonlinear  example, the classical model of propagation of gases in porous media has been
considerer, that is defined by the equation
\begin{equation}
u_t = \Delta\left(u^m\right) =  \text{div}\left(mu^{m-1}\nabla u\right).
\end{equation}

\begin{figure}
\begin{center}

\includegraphics[width=11cm]{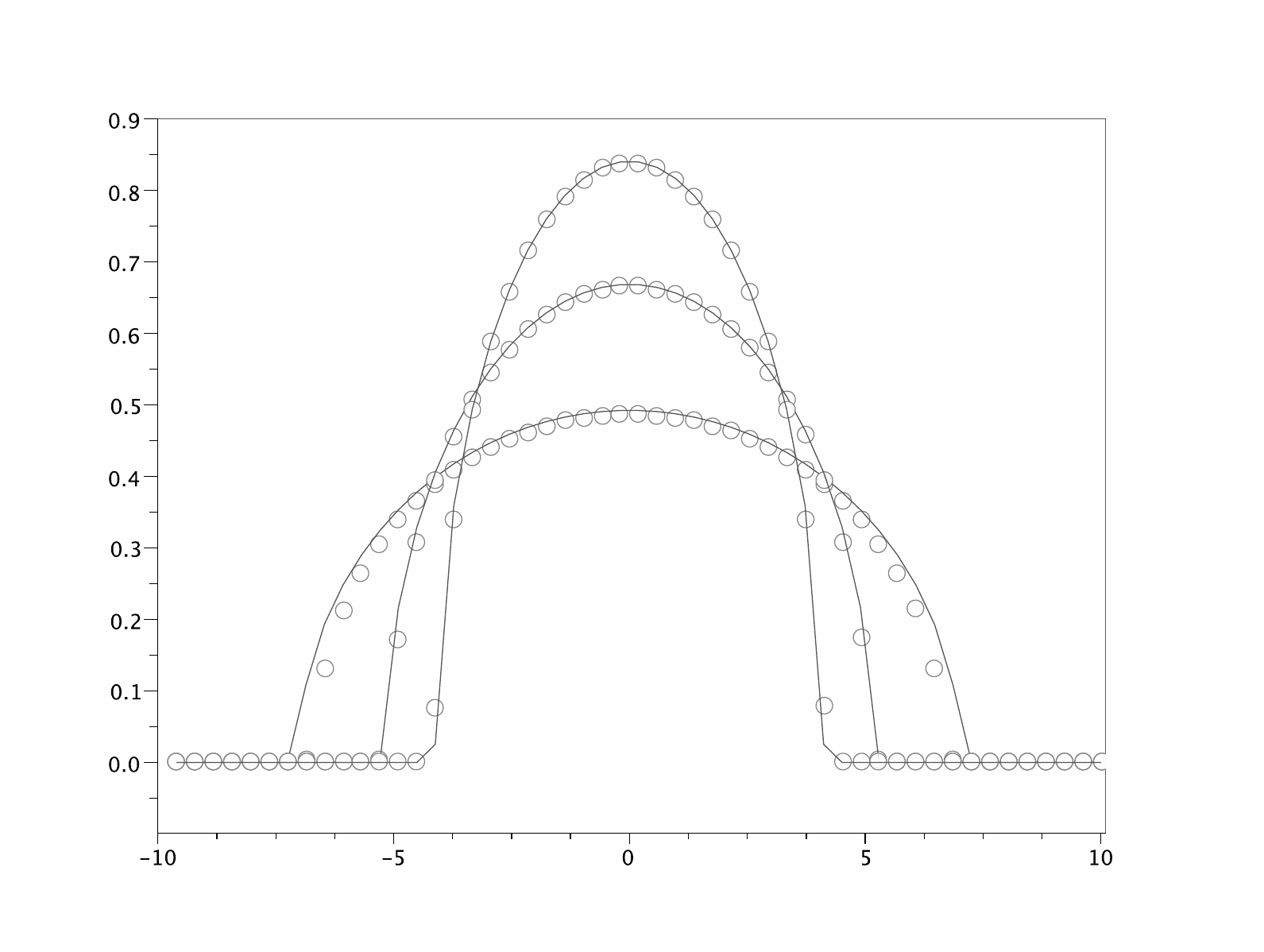}

\caption{Evolution of the exact versus numerical Barenblatt--Pattle solution for T=1,4,16. Numerical solution by SL method (circles), reference solution (continuous line).}
\label{fig:mp1}
\end{center}
\end{figure}

\begin{figure}
\begin{center}

\includegraphics[width=11cm]{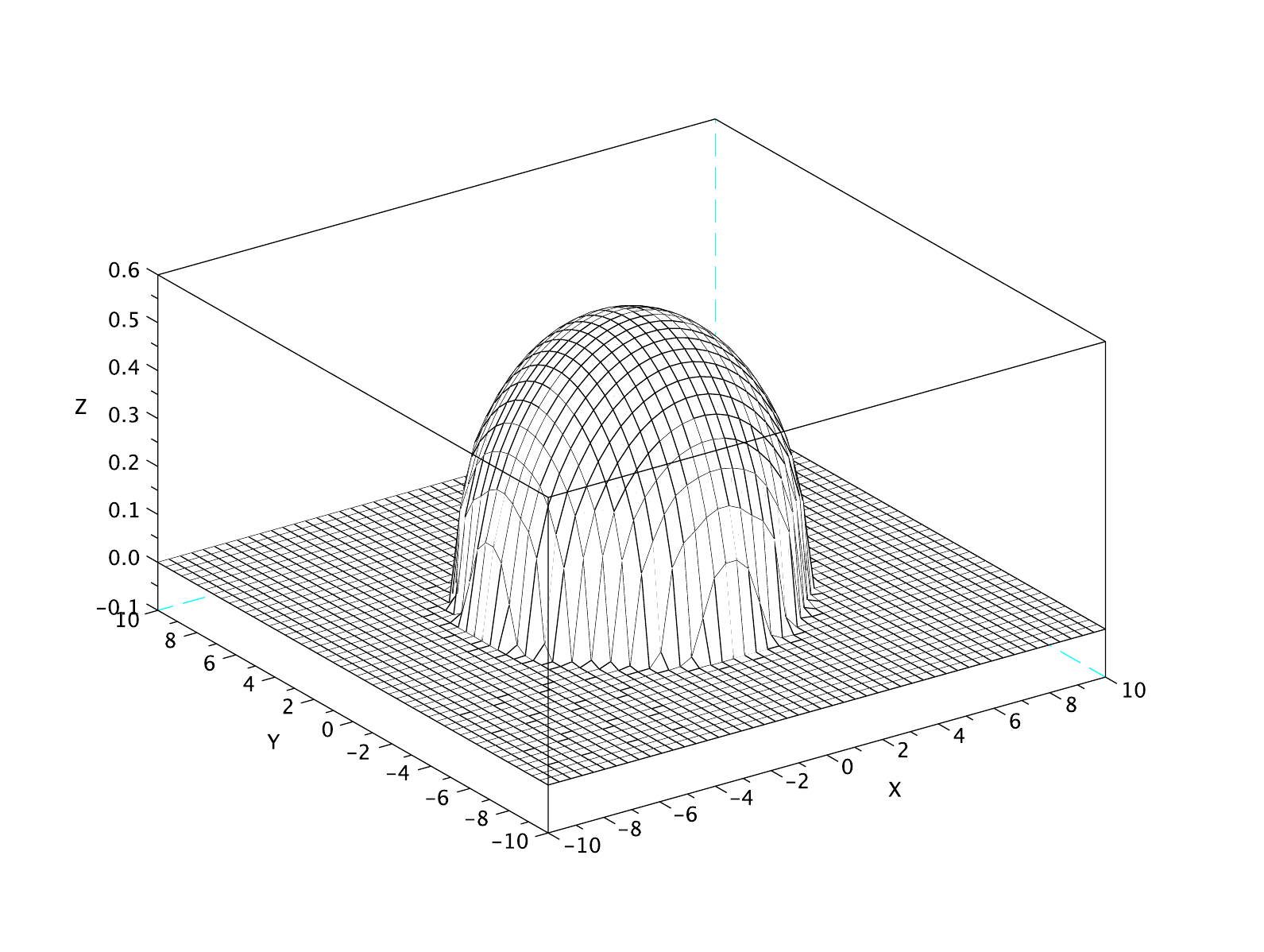}

\caption{Two-dimensional Barenblatt--Pattle solution}
\label{fig:mp2}
\end{center}
\end{figure}

This model   is known  to allow  solutions  that propagate with finite speed. In particular, a family of compactly supported 
self-similar solution has been found independently
 by Barenblatt and Pattle (see, e.g., \cite{barenblatt1952self}), that can be written in the form
\begin{equation}
u(x,t) = (t+t_0)^{-k} \left(A^2 - \frac{k(m-1)|x|^2}{2Md(t+t_0)^{2k/d}}\right)_+^\frac{1}{m-1}
\end{equation}
where $t_0>0$, $A$ is an arbitrary nonzero constant and
\[
k = \frac{d}{d(m-1)+2}.
\]
Figure \ref{fig:mp1} shows the exact and approximate evolution of the Barenblatt--Pattle solution in one space dimension, for $m=3$, $A=1$ and $t_0=1$, computed at $T=1,4,16$ on a mesh composed of 51 nodes, using cubic interpolation and $\Delta t=0.05$. We note that the solution is reasonably accurate even with the rather coarse discretization
employed. Furthermore, the scheme presents no spurious viscosity, since, just like the analytical solution, the numerical solution is also compactly supported. 
Table \ref{table_bp} shows the numerical errors in the 2-norm, at $T=16$, with linear and cubic interpolation under a linear refinement law. The convergence rate is clearly affected by the low regularity (the solution is only H\"older continuous, with exponent $1/2$), but the scheme proves to be quite robust. 
Notice that the computation of the $\delta^\pm_i$ has been carried out iteratively via \eqref{iteraz}  starting from  initial guesses
with large values. Indeed, the iterative computation scheme could get stalled at 
points out of the support of the solution (but in its neighbourhood), due to the possible convergence of \eqref{iteraz} to the fixed point $\delta^\pm_i=0$.
Figure \ref{fig:mp2} shows the result of an analogous computation in two dimensions, at $T=4$, on a mesh composed of $51\times 51$ nodes, using cubic interpolation
and $\Delta t=0.05$.

 \begin{table}[t]
 \centering
\begin{tabular}{||c|c|c|c||} 
\hline
\multicolumn{2}{||c|}{Resolution} & \multicolumn{2}{c||}{$l_2$ relative error} \\
\hline
$N$ & $M$ & $\mathbb{P}_1$ & cubic \\
\hline \hline
50 & 320 & 0.316 & $8.69\cdot 10^{-2}$     \\
\hline
100 & 640 &  0.212 & $4.76\cdot 10^{-2}$     \\
\hline
200 & 1280 &  0.171 & $4.84\cdot 10^{-2}$    \\
\hline
400 & 2560 & 0.135 & $4.49\cdot 10^{-2}$   \\
\hline
800 & 5120 & 0.107 & $3.22\cdot 10^{-2}$   \\
\hline
\end{tabular}
\caption{Relative errors for the Barenblatt--Pattle solution in the 2-norm, linear and cubic interpolation.}
 \label{table_bp}   
 \end{table}


\subsection{Turbulent vertical diffusion}
\label{turb}

As a final nonlinear test, we consider a  simplified model for turbulent vertical diffusion in the atmosphere. Models of this kind
date back as early as \cite{louis:1979} and are widely used in numerical weather prediction models.
The specific version we consider has been employed by several authors (see \cite{girard:1990}, \cite{teixeira:1999}, \cite{benard:2000}) as an example of model
problem closer to the real problems encountered in environmental modelling.
This model is defined by the system of coupled nonlinear diffusion equations
\begin{equation}
\begin{cases}
u_t = (k(u_z,\theta_z)u_z)_z \\
\theta_t = (k(u_z,\theta_z)\theta_z)_z
\end{cases}
\end{equation}
in which $z$ stands for the vertical coordinate, $u$ and $\theta$ denote respectively the wind speed and the potential temperature, and the nonlinear 
diffusivity $k$ has the form
\[
k(u_z,\theta_z) = l^2 |u_z| F(Ri),
\]
with the Richardson number $Ri$ given by
\[
Ri = \frac{g\theta_z}{\theta_0u_z^2}
\]
and $g$ is the acceleration of gravity.
Here, we have taken a mixing length $l=50 $ m, the reference temperature is $\theta_0=273 $ K and $F$ is defined as
\[
F(Ri) = (1+b|Ri|)^\beta
\]
where
\[
\begin{cases}
\beta=-2, \> b=5 & \text{if } Ri>0 \\
\beta=1/2, \> b=20 & \text{if } Ri<0.
\end{cases}
\]
Starting from the same wind speed profile, we show two somewhat typical situations. The first describes a stable configuration in which a first layer with positive temperature gradient is followed by a layer with constant temperature. In the second, unstable configuration, the lower layer has a negative temperature gradient. More precisely, we consider the interval $z\in [0,2000 \ \text{m}]$ with an initial wind speed given, in ms$^{-1}$, by
\[
u_0(z) = 0.2236 z^{1/2},
\]
and boundary conditions $u(0,t)=0$, $u(2000,t)=10$. Initial and boundary conditions for potential temperature, in K, are set respectively as
\[
\begin{cases}
\theta_0(x) = \min(290+0.005 z,295) \\
\theta(0,t) = 290 \ , \ \ \theta(2000,t) = 295
\end{cases}
\]
for the stable case, and
\[
\begin{cases}
\theta_0(x) = \max(300-0.005 z,295) \\
\theta(0,t) = 300 \ , \ \ \theta(2000,t) = 295
\end{cases}
\]
for the unstable case.

\begin{figure}
\begin{center}
\includegraphics[width=6cm]{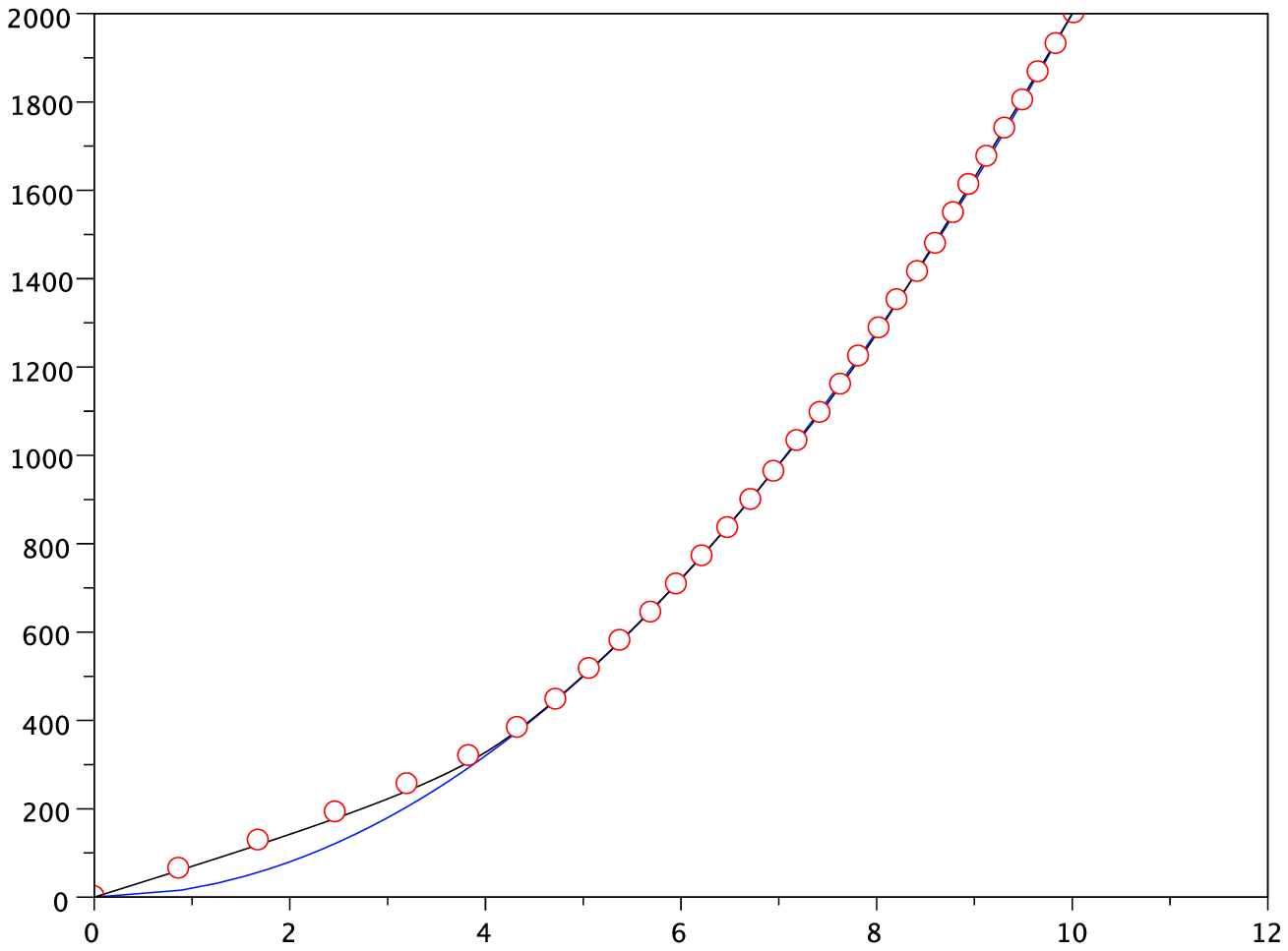}
\includegraphics[width=6cm,height=4.5cm]{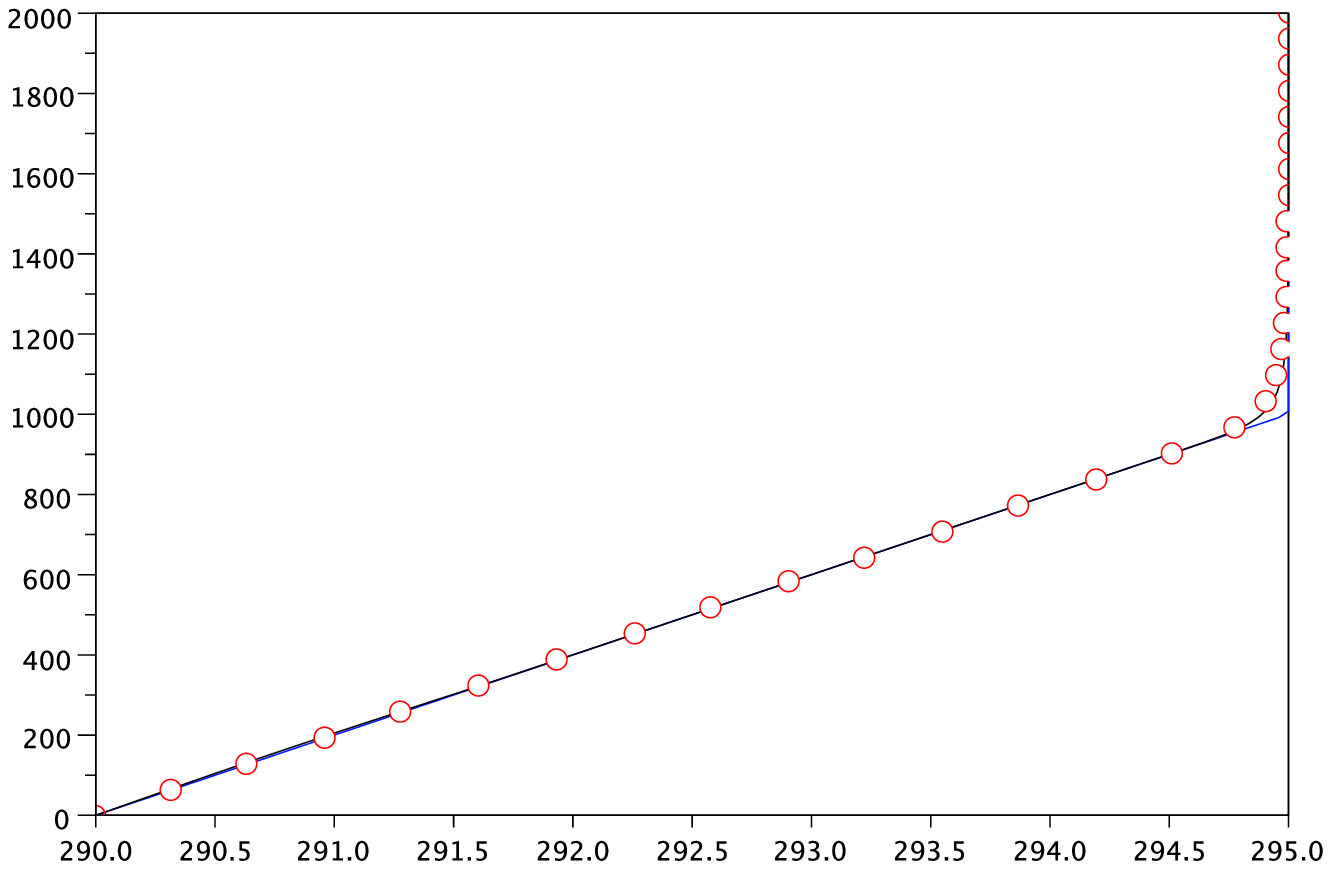}

\includegraphics[width=6cm]{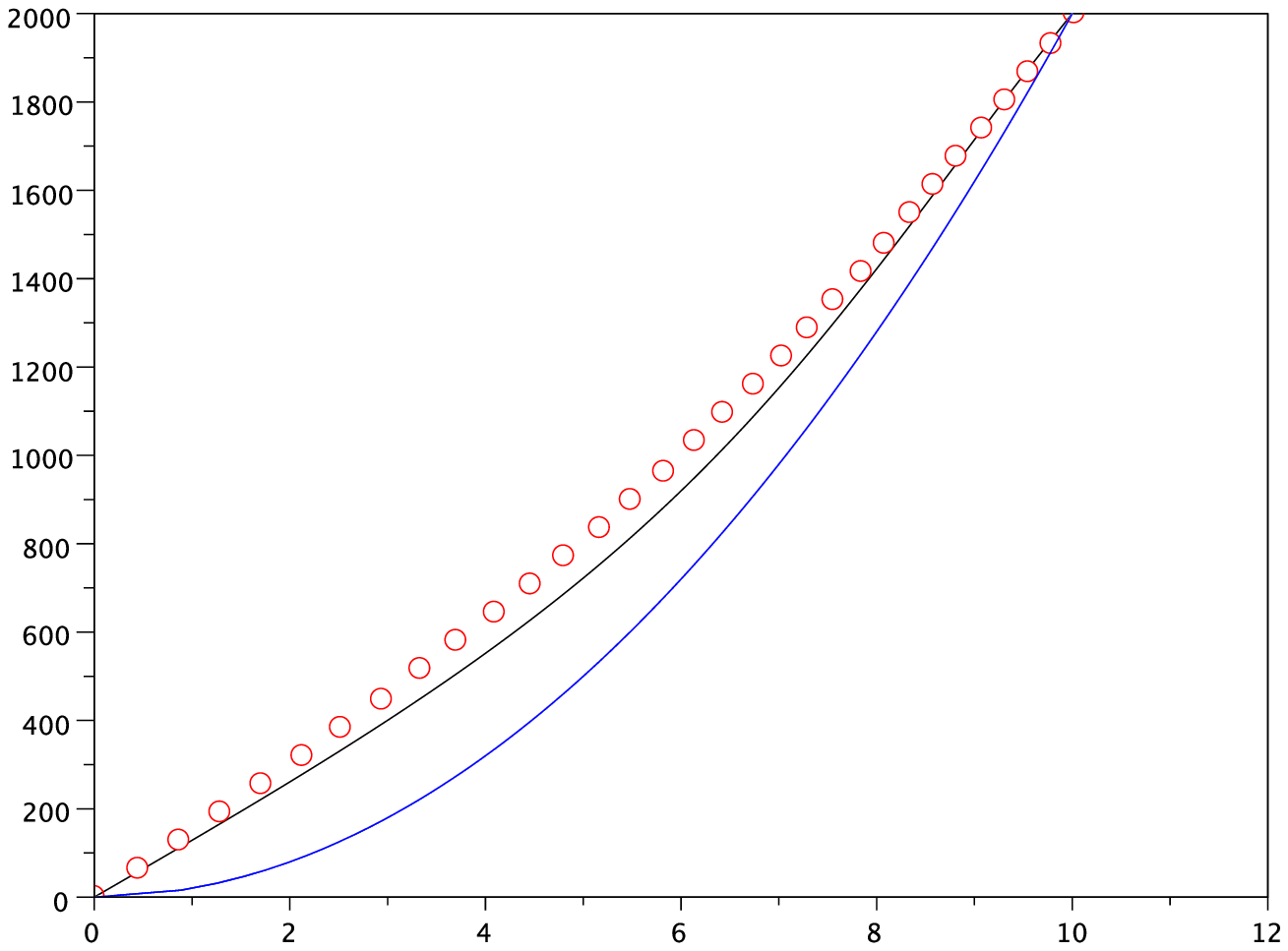}
\includegraphics[width=6cm]{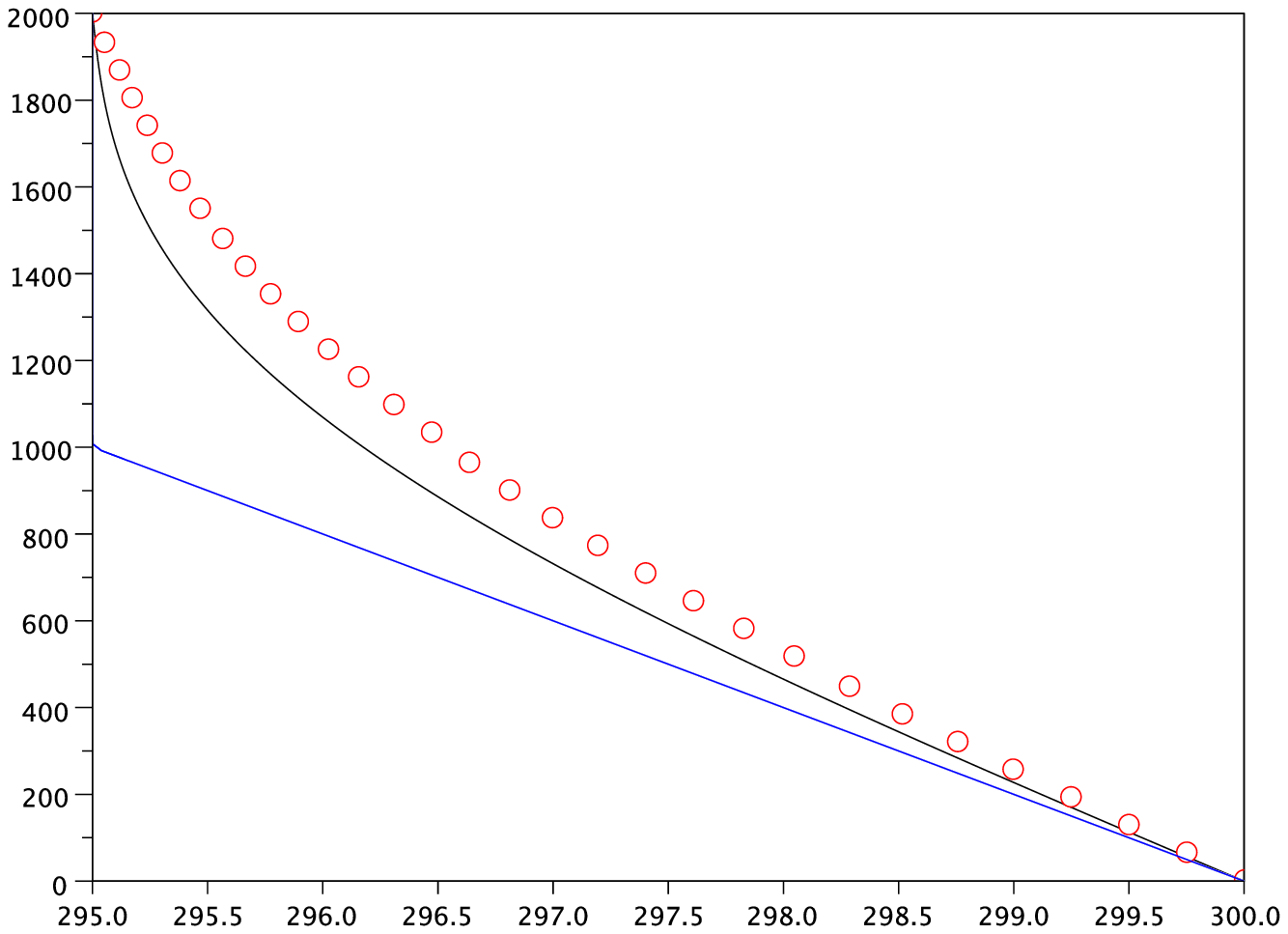}

\caption{Evolution of wind (left) and temperature (right) profiles for the vertical turbulent diffusion model in stable (upper) and unstable (lower) conditions.}
\label{fig:turb}
\end{center}
\end{figure}

Fig. \ref{fig:turb} shows the initial conditions and numerical solutions, computed at the final time of $T=900 $ s, with 32 nodes and $\Delta t=2.5 $ s, with wind speed plots on the left and temperatures on the right. In absence of an exact solution, we compare the numerical result with a reference solution obtained with a four times higher space and time resolution (128 nodes and $\Delta t=0.625 $ s). Note that the qualitative behaviour of the solutions is, as physically expected, a quasi-stationary state in the stable case, and a strong diffusion of the initial condition in the unstable case -- the coarse discretization recovers robustly, although not very accurately, this qualitative behaviour. In this equation, abrupt changes of diffusivity may occur passing from a stable to an unstable layer or vice versa, in which case the fixed point iteration \eqref{iteraz} fails in general to converge. The computation of the $\delta_i^\pm$ has been carried out by a bisection algorithm.


\section{Conclusions and future work}
\label{conclu}

We have extended semi-Lagrangian methods to diffusive problems in divergence form,
providing an advective formulation which is first order consistent in time.
In many computational fluid dynamics applications, the proposed approach has the advantage 
of allowing to handle directly the most common formulations of turbulent diffusivity, without the need to rewrite a 
divergence form diffusive term as an advection--diffusion term with a drift that is due to the variable diffusion coefficient rather than to the physical advection process.
Indeed, in typical anisotropic meshes employed in environmental applications, the diffusion terms
introduce in the differential problem to be solved a remarkable stiffness, which is usually
handled by coupling an Eulerian discretization in space to an implicit time stepping method,
while the proposed approach is explicit and leads easily to an increase of the spatial accuracy 
of the resulting discretization of the second order terms.
The major drawback of the proposed approach seems to be its inherent first order time truncation error.
However,  achieving higher order accuracy with unconditionally stable
and robust schemes is not an easy task (see, e.g., \cite{kalnay:1988}, \cite{teixeira:1999}, \cite{wood:2007}) and the approach we propose fits well in a coherent
SL framework while achieving the same effective accuracy displayed by more standard discretizations. 
A natural development of the approach proposed in this paper is the extension to fully conservative SL schemes
for the advection--diffusion equation and its application to strongly nonlinear diffusion problems such as the Richards
equation for water flow in porous media.

\section*{Acknowledgements}

This research work has been financially supported by the INDAM--GNCS projects
\textit{Metodologie teoriche ed applicazioni avanzate nei metodi Semi-Lagrangiani} and
 \textit{Metodi ad alta risoluzione per problemi evolutivi fortemente nonlineari}, and by Politecnico di Milano.  
We would like to acknowledge useful discussions on the topics analysed in this paper with
O. Bokanowski, E.Sonnendrucker and M.Restelli.

\bibliographystyle{plain}
\bibliography{diff_sl}

\end{document}